\documentclass{amsart}
\usepackage{amsfonts}
\usepackage{amsmath}
\usepackage{amssymb}
\usepackage{amsthm}
\newtheorem{theorem}{Theorem}
\newtheorem{prop}{Proposition}
\newtheorem{lemma}{Lemma}
\newtheorem{rem}{Remark}

\newtheorem{fact}{Fact}
\newtheorem*{conj}{Conjecture}

\begin{document}

\author{Mark Pankov}
\title{Base subsets of polar Grassmannians}
\begin{abstract}
Let $\Delta$ be a thick building of type $\textsf{X}_{n}=\textsf{C}_{n},\textsf{D}_{n}$.
Let also ${\mathcal G}_k$ be the Grassmannian of $k$-dimensional singular subspaces
of the associated polar space $\Pi$ (of rank $n$).
We write ${\mathfrak G}_k$
for the corresponding shadow space of type $\textsf{X}_{n,k}$.
Every bijective transformation of ${\mathcal G}_k$
which maps base subsets to base subsets (the shadows of apartments)
is a collineation of ${\mathfrak G}_k$,
and it is induced by a collineation of
$\Pi$ if $n\ne 4$ or $k\ne 1$.
\end{abstract}

\address{Department of Mathematics and Information Technology,
University of Warmia and Mazury, {\. Z}olnierska 14A, 10-561
Olsztyn, Poland} \email{pankov@matman.uwm.edu.pl}

\maketitle

\section{Introduction}
Let us consider a building of type $\textsf{X}_{n}$ \cite{Tits}
as an incidence geometry $({\mathcal G},*,d)$,
where ${\mathcal G}$ is a set of subspaces, $*$ is an incidence relation,
and $d:{\mathcal G}\to \{0,\dots,n-1\}$ is a dimension function.
Then apartments of the building are certain subgeometries of rank $n$.
Denote by ${\mathcal G}_{k}$ the Grassmannian consisting of all $k$-dimensional subspaces.
The {\it shadow space} ${\mathfrak G}_{k}$ is the partial linear space
whose point set is ${\mathcal G}_{k}$ and lines are defined
by pairs of incident subspaces $S\in {\mathcal G}_{k-1}$, $U\in {\mathcal G}_{k+1}$
if $0<k<n-1$; in the case when $k=0,n-1$,
lines are defined by elements of ${\mathcal G}_1$ and ${\mathcal G}_{n-2}$,
respectively.
The intersection of an apartment with ${\mathcal G}_{k}$
is known as the {\it shadow} of this apartment.

Every apartment preserving (in both directions)
bijective transformation of the chamber set
(the set of maximal flags of our geometry)
can be extended to a bijective transformation of ${\mathcal G}$
preserving the incidence relation \cite{AVM}.

\begin{conj}{\rm
Every bijective transformation of ${\mathcal G}_{k}$
sending the shadows of apartments to the shadows of apartments
is a collineation of ${\mathfrak G}_{k}$.
}\end{conj}

This conjecture holds for the shadow spaces of type $\textsf{A}_{n,k}$ \cite{Pankov2}
and all shadow spaces of symplectic buildings \cite{Pankov3}.
In the present paper we show that this is true
for all buildings of types $\textsf{C}_{n}$ and $\textsf{D}_{n}$
(the buildings associated with polar spaces).

Suppose that our building is of type $\textsf{C}_{n}$ and $\textsf{D}_{n}$.
Let $\Pi$ be the associated polar space (of rank $n$).
We say that $\Pi$ is of type \textsf{C} or \textsf{D}
if the correspondent possibility is realized:
\begin{enumerate}
\item[(\textsf{C})]
every $(n-2)$-dimensional (second maximal) singular subspace
is contained in at least $3$ distinct maximal singular subspaces,
\item[(\textsf{D})]
every $(n-2)$-dimensional singular subspace
is contained in precisely $2$ maximal singular subspaces.
\end{enumerate}
A subset of ${\mathcal G}_{k}$ is said to be {\it inexact}
if it is contained in the shadows of two distinct apartments.
In the case (\textsf{C}),
there are precisely two possibilities for a maximal inexact subset;
as in \cite{Pankov3}, we use maximal inexact subsets of first and second types
to characterize the collinearity relation of ${\mathfrak G}_{k}$
for $k\le n-2$ and $k=n-1$, respectively.
In the case (\textsf{D}),
the first possibility does not realize and we characterize
the collinearity relation in terms of maximal inexact subsets of second type.

\section{Polar Geometry}

\subsection{Partial linear spaces}
Let us consider a pair $\Pi=(P,{\mathcal L})$, where $P$ is a
non-empty set and ${\mathcal L}$ is a family of proper subsets of $P$.
Elements of $P$ and ${\mathcal L}$ will be called {\it points}
and {\it lines}, respectively.
Two or more points are said to be
{\it collinear} if there is a line containing them.
We suppose that $\Pi$ is a {\it partial linear space}:
each line contains at least two points,
for each point there is a line containing it,
and for any two distinct collinear points $p,q\in P$
there is precisely one line containing them,
this line will be denoted by $p\,q$.
We say that $S\subset P$ is a {\it subspace} of $\Pi$ if for any two
distinct collinear points $p,q\in S$ the line $p\, q$ is contained
in $S$.
By this definition, any set consisting of non-collinear
points is a subspace.
A subspace is said to be {\it singular} if any two points of the subspace are collinear
(the empty set and one-point subspaces are singular).

A bijective transformation of $P$
is a {\it collineation} of $\Pi$ if
${\mathcal L}$ is preserved in both directions.

\subsection{Polar spaces}
A {\it polar space of rank} $n$ can be defined as a partial linear space
$\Pi=(P,{\mathcal L})$ satisfying the following axioms:
\begin{enumerate}
\item[(1)] if $p\in P$ and $L\in {\mathcal L}$
then $p$ is collinear with one or all points of $L$,
\item[(2)]
there is no point collinear with all others,
\item[(3)]
every maximal flag of singular subspaces consists of $n$ elements.
\end{enumerate}
It follows from the equivalent axiom system \cite{Tits}
that every singular subspace is contained in a maximal singular subspace and
the restrictions of $\Pi$ to all maximal singular subspaces are
$(n-1)$-dimensional projective spaces
(throughout the paper the dimension is always assumed to be projective).
The complete list of polar spaces of rank greater than $2$
can be found in \cite{Tits}.

The collinearity relation will be denoted by $\perp$:
we write $p\perp q$ if $p$ and $q$ are collinear points,
and $p\not\perp q$ otherwise.
More general, if $X$ and $Y$ are subsets of $P$
then $X\perp Y$ means that every point of $X$ is collinear with all points of $Y$.
The minimal subspace containing $X$ is
called {\it spanned} by $X$ and denoted by $\overline{X}$; in the case when $X\perp X$, it is singular.
For a subset $X\subset P$  we denote by $X^{\perp}$
the set of points collinear with each point of $X$;
it follows from the axiom (1) that $X^{\perp}$ is a subspace.
If $p\perp q$ then
$$\{p,q\}^{\perp\perp}=p\,q^{\perp\perp}=p\,q$$
and we get the following.

\begin{fact}\label{f-1}
Every bijective transformation of $P$ preserving
the collinearity relation {\rm (}in both directions{\rm )}
is a collineation of $\Pi$.
\end{fact}

\begin{rem}{\rm
If $\Pi$ is defined by a sesquilinear or quadratic form
then every collinea\-tion of $\Pi$ can be extended to a collineation
of the corresponding projective space \cite{Chow,Die1},
see also \cite{Die2}.
}\end{rem}

\subsection{Grassmann spaces}
Let $\Pi=(P,{\mathcal L})$ be a polar space of rank $n\ge 3$.
For each $k\in \{0,1,\dots, n-1\}$ we denote by ${\mathcal G}_{k}$
the {\it Grassmannian} consisting of all $k$-dimensional
singular subspaces.
Thus ${\mathcal G}_{0}=P$ and ${\mathcal G}_{1}={\mathcal L}$.

We say that a subset of ${\mathcal G}_{n-1}$ is a {\it line}
if it consists of all maximal singular subspaces containing
a certain $(n-2)$-dimensional singular subspace.
In the case when $1\le k\le n-2$,
a subset of ${\mathcal G}_{k}$ is a {\it line} if
there exist incident $S\in {\mathcal G}_{k-1}$ and $U\in {\mathcal G}_{k+1}$
such that it is the set of all
$k$-dimensional singular subspaces incident with both $S$ and $U$.
The family of all lines in ${\mathcal G}_{k}$ will be denoted by ${\mathcal L}_{k}$.
The pair
$${\mathfrak G}_{k}:=({\mathcal G}_{k},{\mathcal L}_{k})$$
is a partial linear space
for each $k\in\{1,\dots,n-1\}$.
In what follows we will suppose that ${\mathfrak G}_{0}=\Pi$.

Two distinct maximal singular subspaces are collinear points of
${\mathfrak G}_{n-1}$ if and only if
their intersection belongs to ${\mathcal G}_{n-2}$.
In the case when $k\le n-2$, two distinct
$k$-dimensional singular subspaces $S$ and $U$ are collinear points of
${\mathfrak G}_{k}$ if and only if $S\perp U$
and the subspace spanned by them belongs to ${\mathcal G}_{k+1}$.
We say that two elements of ${\mathcal G}_{k}$ ($k\le n-2$)
are {\it weak-adjacent} if their intersection is an element of ${\mathcal G}_{k-1}$;
it is trivial that any two distinct collinear points of ${\mathfrak G}_{k}$
are weak-adjacent, the converse fails.

Every collineation of $\Pi$ induces a collineation of ${\mathfrak G}_{k}$.
Conversely,

\begin{fact}\label{f-2}
Every bijective transformation of ${\mathcal G}_{k}$ $(1\le k\le n-1)$
preserving the collinearity relation {\rm (}in both directions{\rm )}
is a collineation of ${\mathfrak G}_{k}$;
moreover, it is induced by a collineation of $\Pi$ if
$n\ne 4$ or $k\ne 1$.
Every bijective transformation of ${\mathcal G}_{k}$ $(1\le k\le n-2)$
preserving the weak-adjacency relation {\rm (}in both directions{\rm )}
is  the collineation of ${\mathfrak G}_{k}$ induced by a collineation of $\Pi$.
\end{fact}

\begin{rem}{\rm
It was proved in \cite{Chow,Die1} for $k=n-1$ (see also \cite{Die2}),
and we refer \cite{PPZ} for the general case.
}\end{rem}

Now suppose that $n\ge 4$ and the case (\textsf{D})
is realized.
Then $\Pi$ is the polar space of a non-degenerate quadratic from $q$
defined on a certain $2n$-dimensional vector space
(if the characteristic of the field is not equal to $2$ then
the associated bilinear form is symmetric and non-degenerate,
and the corresponding polar space coincides with $\Pi$).
In this case, each line of ${\mathfrak G}_{n-1}$ consists of
precisely $2$ points.
Let ${\mathcal O}_{+}$ and ${\mathcal O}_{-}$
be the orbits of the action of the orthogonal group ${\rm O}^{+}(q)$ on
${\mathcal G}_{n-1}$.
Then
$$n-\dim (S\cap U)\;\mbox{ is }\;
\begin{cases}
\mbox{ odd }&  \mbox{ if } S,U\in {\mathcal O}_{\delta},\;\delta=+,-\\
\mbox{ even }& \mbox{ if } S\in {\mathcal O}_{+},\;U\in {\mathcal O}_{-}.
\end{cases}
$$
A subset of ${\mathcal O}_{\delta}$ is called a {\it line}
if it consists of all elements of ${\mathcal O}_{\delta}$
containing a certain $(n-3)$-dimensional singular subspace.
We denote by ${\mathcal L}_{\delta}$
the family of all lines and get the partial linear space
$${\mathfrak O}_{\delta}=({\mathcal O}_{\delta},{\mathcal L}_{\delta}).$$
Two distinct elements of ${\mathcal O}_{\delta}$ are collinear points of
${\mathfrak O}_{\delta}$ if and only if their intersection is
an element of ${\mathcal G}_{n-3}$.
Every collineation of $\Pi$ induces collineations of ${\mathfrak O}_{+}$
and ${\mathfrak O}_{-}$.

\begin{fact}\label{f-3}\cite{Chow,Die1}
Every bijective transformation of ${\mathcal O}_{\delta}$ $(\delta =+,-)$
preserving the collinearity relation {\rm (}in both directions{\rm )}
is a collineation of ${\mathfrak O}_{\delta}$ induced by a collinea\-tion of $\Pi$.
\end{fact}

\begin{rem}{\rm
Singular subspaces of the partial linear spaces
considered above are well known,
some information concerning non-singular subspaces
can be found in \cite{CKS}.
}\end{rem}

\subsection{Base subsets}
We say that $B=\{p_{1},\dots,p_{2n}\}$ is a {\it base} of $\Pi$
if for each $i\in \{1,\dots,2n\}$
there exists unique $\sigma(i)\in \{1,\dots,2n\}$ such that
$$p_{i}\not\perp p_{\sigma(i)}.$$
In this case, the set of all $k$-dimensional singular subspaces spanned by points of $B$
is said to be the {\it base subset} of ${\mathcal G}_{k}$ associated with (defined by) $B$;
it coincides with $B$ if $k=0$.

Every base subset of ${\mathcal G}_{k}$ consists of precisely
$$2^{k+1}\binom{n}{k+1}$$
elements (Proposition 1 in \cite{Pankov3}).

\begin{prop}\label{prop1}
For any two $k$-dimensional singular subspaces of $\Pi$ there is a
base subset of ${\mathcal G}_{k}$ containing them.
\end{prop}

\begin{proof}
Easy verification.
\end{proof}

We define {\it base subsets} of ${\mathcal O}_{\delta}$
as the intersections of ${\mathcal O}_{\delta}$ with base subsets of ${\mathcal G}_{n-1}$.
Every base subset of ${\mathcal O}_{\delta}$ consists of $2^{n-1}$ elements
and it follows from Proposition \ref{prop1} that for any two elements of
${\mathcal O}_{\delta}$ there is a base subset of ${\mathcal O}_{\delta}$ containing
them.

\subsection{Buildings associated with polar spaces}
In the case (\textsf{C}),
the incidence geometry of singular subspaces is a thick building of type
$\textsf{C}_{n}$,
where each apartment is the subgeometry consisting
of all singular subspaces spanned by points of a certain base of $\Pi$;
the shadow spaces are ${\mathfrak G}_{k}$ ($0\le k\le n-1$),
and the shadows of apartments are base subsets.

Suppose that $n\ge 4$ and $\Pi$ is of type \textsf{D}.
The {\it oriflamme} incidence geometry consists of
all singular subspaces of dimension distinct from $n-2$;
two singular subspaces $S$ and $U$ are incident
if $S\subset U$, or $U\subset S$, or the dimension of
$S\cap U$ is equal to $n-2$.
This is a thick building of type $\textsf{D}_{n}$;
as above, apartments are the subgeometries defined by bases of $\Pi$.
The shadow spaces are ${\mathfrak G}_{k}$ ($0\le k\le n-3$)
and ${\mathfrak O}_{\delta}$ ($\delta=+,-$),
the shadows of apartments are base subsets.

\section{Results}

\begin{theorem}\label{theorem1}
Let $f$ be a bijective transformation of ${\mathcal G}_{k}$ $(0\le k\le n-1)$
which maps base subsets to base subsets
\footnote{We do not require that $f$ preserves
the class of base subsets in both directions.}.
Then $f$ is a collineation of ${\mathfrak G}_{k}$.
If $k\ge 1$ and $\Pi$ is of type {\rm \textsf{C}}
then this collineation is induced by a collineation of $\Pi$.
In the case when $k\ge 1$ and $\Pi$ is of type {\rm \textsf{D}},
it is induced by a collineation of $\Pi$ if $n\ne 4$ or $k\ne 1$.
\end{theorem}

\begin{theorem}\label{theorem5}
If $n\ge 4$ and $\Pi$ is of type {\rm \textsf{D}}
then every bijective transformation of ${\mathcal O}_{\delta}$
$(\delta=+,-)$
which maps base subsets to base subsets
\footnote{As in Theorem 1, we do not require that
our transformation preserves the class of base subsets in both directions.}
is the colleniation of ${\mathfrak O}_{\delta}$
induced by a collineation of $\Pi$.
\end{theorem}

\section{Inexact and complement subsets}

\subsection{Inexact subsets}
Let $B=\{p_{1},\dots,p_{2n}\}$ be a base of $\Pi$
and ${\mathcal B}$ be the associated base subset of ${\mathcal G}_{k}$.
Recall that ${\mathcal B}$ consists of all $k$-dimensional singular subspaces
spanned by $p_{i_{1}},\dots, p_{i_{k+1}}$,
where
$$\{i_{1},\dots,i_{k+1}\}\cap \{\sigma(i_{1}),\dots,\sigma(i_{k+1})\}=\emptyset.$$
If $k=n-1$ then every element of
${\mathcal B}$ contains precisely one of the points $p_{i}$ or $p_{\sigma(i)}$
for each $i$.

We write ${\mathcal B}(+i)$ and ${\mathcal B}(-i)$ for the sets of all
elements of ${\mathcal B}$ which  contain $p_{i}$ or do not contain $p_{i}$, respectively.
For any $i_{1},\dots, i_{s}$ and $j_{1},\dots, j_{u}$
belonging to $\{1,\dots,2n\}$
we define
$${\mathcal B}(+i_{1},\dots, +i_{s},-j_{1},\dots, -j_{u}):=
{\mathcal B}(+i_{1})\cap\dots\cap {\mathcal B}(+i_{s})\cap
{\mathcal B}(-j_{1})\cap \dots\cap{\mathcal B}(-j_{u}).$$
It is trivial that
$${\mathcal B}(+i)={\mathcal B}(+i,-\sigma(i))$$
and in the case when $k=n-1$ we have
$${\mathcal B}(+i)={\mathcal B}(-\sigma(i)).$$

Let ${\mathcal R}\subset {\mathcal B}$.
We say that ${\mathcal R}$ is {\it exact} if
there is only one base subset of ${\mathcal G}_{k}$ containing ${\mathcal R}$;
otherwise, ${\mathcal R}$ will be called {\it inexact}.
If ${\mathcal R}\cap {\mathcal B}(+i)$ is not empty then we define
$S_{i}({\mathcal R})$ as the intersection of all subspaces
belonging to ${\mathcal R}$ and containing $p_{i}$,
and we define $S_{i}({\mathcal R}):=\emptyset$ if
the intersection of ${\mathcal R}$ and ${\mathcal B}(+i)$ is empty.
If
$$S_{i}({\mathcal R})=p_{i}$$
for all $i$ then ${\mathcal R}$ is exact;
the converse fails.

\begin{prop}\label{prop2-1}
If $k=n-1$ then ${\mathcal B}(-i)$ is inexact, but this inexact subset is not maximal.
In the case when $k\in \{0,\dots,n-2\}$,
the following assertions are fulfilled:
\begin{enumerate}
\item[$\bullet$]
if $\Pi$ is of type {\rm \textsf{C}} then
${\mathcal B}(-i)$ is a maximal inexact subset,
\item[$\bullet$]
${\mathcal B}(-i)$ is exact
if $\Pi$ is of type {\rm \textsf{D}}.
\end{enumerate}
\end{prop}

\begin{proof}
If $k=n-1$ then
for each $U\in {\mathcal B}\setminus {\mathcal B}(-i)$
the subset
$${\mathcal B}(-i)\cup\{U\}$$
is inexact
(see proof of Proposition 3 in \cite{Pankov3} for details);
this means that ${\mathcal B}(-i)$ is a non-maximal inexact subset.

Let $k\le n-2$.
Then
$$S_{j}({\mathcal B}(-i))=p_{j}\;\mbox{ for }\;j\ne i$$
(proof of Proposition 3 in \cite{Pankov3}).
Therefore, if ${\mathcal B}(-i)$ is contained
in the base subset of ${\mathcal G}_{k}$ associated with
a base $B'$ of $\Pi$ then
\begin{equation}\label{equation2-1}
B'=(B\setminus\{p_{i}\})\cup \{p\}.
\end{equation}
It is clear that $p$ is collinear with all points of $B\setminus\{p_{i},p_{\sigma(i)}\}$
and non-collinear with $p_{\sigma(i)}$.
If $p\ne p_{i}$
then $p\not\perp p_{i}$ and each $(n-2)$-dimensional singular subspace spanned by points
of $B\setminus\{p_{i},p_{\sigma(i)}\}$
is contained in $3$ distinct maximal singular subspaces which contradicts
(\textsf{D}).

Suppose that $\Pi$ is of type \textsf{C}.
Let $S$ and $U$ be non-intersecting $(n-2)$-dimensional
singular subspaces spanned by points of $B\setminus\{p_{i},p_{\sigma(i)}\}$.
By (\textsf{C}), there exists a point $p\in P$
non-collinear with $p_{i},p_{\sigma(i)}$
and collinear with all points of $S$.
Since any two flags are contained in a certain apartment,
there is a point $q\in P$ such that
$$(B\setminus\{p_{i},p_{\sigma(i)}\})\cup\{p,q\}$$
is a base of $\Pi$ whose points span $U$ and $\overline{S\cup\{p\}}$.
This means that $p$ is collinear with all points of $B\setminus\{p_{i},p_{\sigma(i)}\}$.
The base subset of ${\mathcal G}_{k}$ associated with the base
\eqref{equation2-1} contains ${\mathcal B}(-i)$.
So it is an inexact subset.
As in \cite{Pankov3} (proof Proposition 3),
we show that this inexact subset is maximal.
\end{proof}

\begin{prop}\label{prop2-2}
If $j\ne i,\sigma(i)$ then
$${\mathcal R}_{ij}:={\mathcal B}(+i,+j)\cup{\mathcal B}(+\sigma(i),+\sigma(j))
\cup{\mathcal B}(-i,-\sigma(j))$$
is inexact;
moreover, it is a maximal inexact subset except the case when $k=0$ and
$\Pi$ is of  type {\rm \textsf{C}}.
\end{prop}

\begin{proof}
As in \cite{Pankov3} (proof of Proposition 4),
we establish that ${\mathcal R}_{ij}$ is a maximal inexact subset if $k\ge 1$.
If $k=0$ then
$${\mathcal R}_{ij}=B\setminus\{p_{i},p_{\sigma(j)}\}.$$
We choose a point $p'_{i}$ on the line $p_{i}p_{j}$ distinct from $p_{i}$ and $p_{j}$,
the line  $p_{\sigma(j)}p_{\sigma(i)}$ contains
the unique point $p'_{\sigma(j)}$ collinear with $p'_{i}$.
Then
$$(B\setminus\{p_{i},p_{\sigma(j)}\})\cup \{p'_{i},p'_{\sigma(j)}\}$$
is a base of $\Pi$
and ${\mathcal R}_{ij}$ is inexact.
If $\Pi$ is of  type {\rm \textsf{C}}
then ${\mathcal B}(-i)=B\setminus\{p_{i}\}$
is an inexact subset containing ${\mathcal R}_{ij}$,
hence the inexact subset ${\mathcal R}_{ij}$ is not maximal.
\end{proof}

A direct verification shows that
$${\mathcal R}_{ij}={\mathcal B}(+i,+j)\cup{\mathcal B}(-i)$$
if $k=n-1$.

As in \cite{Pankov3}, the inexact subsets considered in
Propositions \ref{prop2-1} and \ref{prop2-2}
will be called of {\it first} and {\it second} type, respectively.

\begin{prop}\label{prop2-3}
Every maximal inexact subset is of first or second type.
In particular, each maximal inexact subset is of second type
if $k=n-1$ or $\Pi$ is of type {\rm \textsf{D}},
and each maximal inexact subset is of first type
if $k=0$ and $\Pi$ is of type {\rm \textsf{C}}.
\end{prop}

\begin{proof}
The case $k=0$ is trivial.
Let $k\ge 1$ and ${\mathcal R}$ be a maximal inexect subset of ${\mathcal B}$.
First, we consider the case when all $S_{i}({\mathcal R})$  are non-empty.
Denote by $I$ the set of all $i$ such that
the dimension of $S_{i}({\mathcal R})$ is non-zero.
Since ${\mathcal R}$ is inexact, $I$ is non-empty.
Suppose that for certain $l\in I$
the subspace $S_{l}({\mathcal R})$ is spanned by $p_{l},p_{j_{1}},\dots, p_{j_{u}}$
and
$$M_{1}:=S_{\sigma(j_{1})}({\mathcal R}),\dots,
M_{u}:=S_{\sigma(j_{u})}({\mathcal R})$$
do not contain $p_{\sigma(l)}$.
Then $p_{l}$ belongs to $M^{\perp}_{1},\dots,M^{\perp}_{u}$;
on the other hand,
$$p_{j_{1}}\not\in M^{\perp}_{1},\dots,p_{j_{u}}\not\in M^{\perp}_{u}$$
and we have
$$M^{\perp}_{1}\cap\dots\cap M^{\perp}_{u}\cap S_{l}({\mathcal R})=p_{l}.$$
If this holds for all $l\in I$
then our subset is exact.
Therefore, there exist $i\in I$ and $j\ne i,\sigma(i)$
such that
$$p_{j}\in S_{i}({\mathcal R})\;\mbox{ and }\;p_{\sigma(i)}\in S_{\sigma(j)}({\mathcal R}).$$
Then ${\mathcal R}$ is contained ${\mathcal R}_{ij}$.
We have ${\mathcal R}={\mathcal R}_{ij}$, since the inexact subset ${\mathcal R}$
is maximal.

Suppose that $S_{i}({\mathcal R})=\emptyset$ for certain $i$.
Then ${\mathcal R}$ is contained in ${\mathcal B}(-i)$
and we get a maximal inexact subset of first type if
$\Pi$ is of type \textsf{C}.
Otherwise, one of the following possibilities is realized:
\begin{enumerate}
\item[(1)]
$S_{j}({\mathcal R})=\emptyset$ for certain $j\ne i,\sigma(i)$.
Then
${\mathcal R}\subset {\mathcal B}(-j,-i)$ is a proper subset of ${\mathcal R}_{j\,\sigma(i)}$.
\item[(2)]
All $S_{l}({\mathcal R})$, $l\ne i,\sigma(i)$ are non-empty
and there exists $j\ne i,\sigma(i)$ such that $S_{j}({\mathcal R})$
contains $p_{\sigma(i)}$.
Then
$${\mathcal R}\subset{\mathcal B}(+j,+\sigma(i))\cup{\mathcal B}(-j,-i);$$
and as in the previous case,
${\mathcal R}$ is a proper subset of ${\mathcal R}_{j\,\sigma(i)}$.
\item[(3)]
Each $S_{j}({\mathcal R})$, $j\ne i,\sigma(i)$ is non-empty
and does not contain $p_{\sigma(i)}$.
As in the case when $S_{j}({\mathcal R})\ne \emptyset$
for all $j\in \{1,\dots,2n\}$,
we show that
${\mathcal R}$ is contained in certain ${\mathcal R}_{l m}$
with $l,m\ne i,\sigma(i)$.
It is clear that
$${\mathcal R}\subset {\mathcal B}(-i)\cap {\mathcal R}_{l m}$$
is a proper subset of ${\mathcal R}_{l m}$.
\end{enumerate}
In each of the cases considered above, ${\mathcal R}$ is a proper subset
of a maximal inexact subset;
this contradicts the fact that our inexact subset is maximal.
\end{proof}

\subsection{Complement subsets}
Let ${\mathcal B}$ be as in the previous subsection.
We say that ${\mathcal R}\subset {\mathcal B}$ is a {\it complement subset}
if ${\mathcal B}\setminus {\mathcal R}$
is a maximal inexact subset.
A complement subset is said to be of {\it first} or {\it second} type
if the corresponding maximal inexact subset is of first or second type,
respectively.
The complement subsets for the maximal inexact subsets
from Propositions \ref{prop2-1} and \ref{prop2-2} are
$${\mathcal B}(+i)\;\mbox{ and }\;
{\mathcal C}_{ij}:={\mathcal B}(+i,-j)\cup{\mathcal B}(+\sigma(j),-\sigma(i)).$$
If $k=n-1$ then the second subset coincides with
$${\mathcal B}(+i,+\sigma(j))={\mathcal B}(+i,+\sigma(j),-j,-\sigma(i)).$$
In the case when $k\ge 1$,
each complement subset is of second type if $k=n-1$ or $\Pi$ is of type {\rm \textsf{D}}.
If $k=0$ and $\Pi$ is of type {\rm \textsf{D}}
then all complement subsets are of second type (pairs of collinear points).
The case when $k=0$ and $\Pi$ is of type {\rm \textsf{C}} is trivial:
each complement subset consists of a single point.

\begin{lemma}\label{lemma2-1}
Suppose that $1\le k\le n-2$ and $\Pi$ is of type {\rm \textsf{C}}.
Let ${\mathcal R}$ be a complement subset of ${\mathcal B}$.
If ${\mathcal R}$ is of first type then there are precisely $4n-3$ distinct complement
subsets of ${\mathcal B}$ which do not intersect ${\mathcal R}$.
If ${\mathcal R}$ is
of second type then there are precisely $4$ distinct complement
subsets of ${\mathcal B}$ which do not intersect ${\mathcal R}$.
\end{lemma}

\begin{proof}
See proof of Lemma 5 in \cite{Pankov3}.
\end{proof}

\begin{lemma}\label{lemma2-2}
Suppose that $k=n-1$.
Then $S,U\in {\mathcal B}$ are collinear points of ${\mathfrak G}_{k}$
if and only if there are precisely
$$\binom{n-1}{2}$$
distinct complement subsets of ${\mathcal B}$ containing $S$ and $U$.
\end{lemma}

\begin{proof}
Let $\dim(S\cap U)=m$.
The complement subset ${\mathcal B}(+i,+j)$ contains $S$ and $U$ if and only if
the line $p_{i}p_{j}$ is contained $S\cap U$. Thus
there are
$$\binom{m+1}{2}$$
distinct complement subsets of ${\mathcal B}$ containing $S$ and $U$.
\end{proof}

\begin{lemma}\label{lemma-new}
Let $1\le k\le n-2$ and $\Pi$ be of type {\rm \textsf{C}}.
If ${\mathcal R}$ is the intersection of $2n-k-2$ distinct complement subsets
of first type of ${\mathcal B}$ then one of the following possibilities is realized:
\begin{enumerate}
\item[$\bullet$] ${\mathcal R}$ consists of $k+2$ mutually collinear
points of ${\mathfrak G}_{k}$,
\item[$\bullet$] ${\mathcal R}$ consists of $2$ weak-adjacent elements of
${\mathcal G}_{k}$ which are non-collinear.
\end{enumerate}
\end{lemma}

\begin{proof}
Direct verification.
\end{proof}

Let $c(S,U)$ be the number of compliment subsets of {\it second} type of ${\mathcal B}$
containing both $S,U\in{\mathcal B}$.

\begin{lemma}\label{lemma2-4}
Suppose that $1\le k\le n-2$.
Let $S$ and $U$ be distinct elements of ${\mathcal B}$
and $\dim (S\cap U)=m$. Then
$$c(S,U)\le (m+1)(2n-2k+m-2)+(k-m)^{2};$$
moreover, if $k=n-2$ then
$$c(S,U)\le (m+1)(m+2)+1.$$
If $S$ and $U$ are collinear points of ${\mathfrak G}_{k}$ then
$$c(S,U)=k(2n-k-3)+1.$$
In the case when $S$ and $U$ are weak adjacent and non-collinear,
we have
$$c(S,U)=k(2n-k-3).$$
\end{lemma}

\begin{proof}
If the complement subset ${\mathcal C}_{ij}$
contains $S$ and $U$ then one of the following possibilities is realized:
\begin{enumerate}
\item[(A)]
at least one of the subsets ${\mathcal B}(+i,-j)$ or ${\mathcal B}(+\sigma(j),-\sigma(i))$
contains both $S,U$;
\item[(B)]
each of the subsets ${\mathcal B}(+i,-j)$ and ${\mathcal B}(+\sigma(j),-\sigma(i))$
contains precisely one of our subspaces.
\end{enumerate}

The case (A).
If $S$ and $U$ both are contained in ${\mathcal B}(+i,-j)$ then
$$p_{i}\in S\cap U\;\mbox{ and }\;p_{j}\not\in S\cup U$$
Since $S\cup U$ contains precisely $2k-m+1$ points of $B$
and $j\ne \sigma(i)$, we have $2n-2k+m-2$
distinct possibilities for $j$ if $i$ is fixed.
Thus there are precisely
$$(m+1)(2n-2k+m-2)$$
complement subsets satisfying (A),
it must be taken into account that
${\mathcal C}_{ij}$ coincides with ${\mathcal C}_{\sigma(j)\sigma(i)}$.

The case (B).
If
$$S\in {\mathcal B}(+i,-j)\setminus {\mathcal B}(+\sigma(j),-\sigma(i))
\;\mbox{ and }\;
U\in{\mathcal B}(+\sigma(j),-\sigma(i))\setminus {\mathcal B}(+i,-j)$$
then
$$p_{i}\in S\setminus U \;\mbox{ and }\;p_{\sigma(j)}\in U\setminus S$$
(since $p_{\sigma(j)}\in U$, we have $p_{j}\not\in U$;
then $U\not\in {\mathcal B}(+i,-j)$
guarantees that $p_{i}\not\in U$; similarly, we show that $p_{\sigma(j)}\not\in S$).
Hence there are at most
$$(k-m)^{2}$$
distinct complement subsets satisfying (B)
(as above, we take into account that ${\mathcal C}_{ij}$ coincides with ${\mathcal C}_{\sigma(j)\sigma(i)}$).
Note that
$$p_{i}\perp U\;\mbox{ and }\;p_{\sigma(j)}\perp S$$
(since $p_{\sigma(i)}\not\in U$ and $p_{j}\not\in S$).
In the case when $k=n-2$, there are at most one point of
$$B\cap(S\setminus U)$$
collinear with all points of $U$ and at most one point of
$$B\cap(U\setminus S)$$
collinear with all points of $S$;
hence there is at most one complement subset satisfying (B).

Now suppose that $m=k-1$.
It was shown above that there are precisely
$k(2n-k-3)$
complement subsets of kind (A).
Let
$$B\cap(S\setminus U)=\{p_{u}\}\;\mbox{ and }\;B\cap(U\setminus S)=\{p_{v}\}.$$
If $S$ and $U$ are collinear points of ${\mathfrak G}_{k}$ then
$v\ne \sigma(u)$ and the complement subset
${\mathcal C}_{u \sigma(v)}={\mathcal C}_{v \sigma(u)}$
contains $S$ and $U$.
Otherwise,
we have $v=\sigma(u)$
and there are no complement subsets of kind (B).
\end{proof}

Let us define
$${\rm m}_{c}:=\max{\{\;c(S,U)\;:\;S,U\in {\mathcal B},\;S\ne U\;\}}.$$

\begin{lemma}\label{lemma2-5}
Distinct $S,U\in {\mathcal B}$ are collinear points of ${\mathfrak G}_{k}$ if and only if
$c(S,U)={\rm m_c}$.
\end{lemma}

\begin{proof}
Since
$$g(x):=(x+1)(2n-2k+x-2)+(k-x)^{2},\;\;x\in {\mathbb R}$$
is a parabola and the coefficient of $x^2$ is positive,
$$\max{\{g(-1),g(0),\dots,g(k-1)\}}=\max{\{g(-1),g(k-1)\}}.$$
Suppose that $k\le n-3$.
Then
$$g(-1)<g(k-1)$$
($g(k-1)=k(2n-k-2)+1\ge k(k+3)+1>(k+1)^{2}=g(-1)$)
and Lemma \ref{lemma2-4} guarantees that
$c(S,U)={\rm m_c}$ if and only if
$S$ and $U$ are collinear points of ${\mathfrak G}_{k}$.
If $k=n-2$ then
$$k(2n-k-3)+1=k(k+1)+1>(m+1)(m+2)+1\;\mbox{ for }\; m\le k-2$$
and Lemma \ref{lemma2-4} gives the claim.
\end{proof}

\subsection{Proof of Theorem 1}
Let $f$ be a bijective transformation of ${\mathcal G}_{k}$
sending base subsets to base subsets.

First we consider the case when $k>0$.
Let $S,U\in {\mathcal G}_{k}$ and ${\mathcal B}$
be a base subset of ${\mathcal G}_{k}$ containing $S$ and $U$.
It is clear that $f$ transfers an inexact subset of ${\mathcal B}$
to inexact subset of $f({\mathcal B})$.
Since the base subsets ${\mathcal B}$  and $f({\mathcal B})$
have the same number of inexact subsets,
every inexact subset of $f({\mathcal B})$ is
the $f$-image of a certain inexact subset of ${\mathcal B}$.
This implies that an inexact subset of
${\mathcal B}$ is maximal if and only if
its $f$-image is a maximal inexact subset of $f({\mathcal B})$.
Thus $f$ and $f^{-1}$ map complement subsets to complement subsets.
Lemma \ref{lemma2-1} guarantees that the type of a complement subset is preserved
\footnote{This is trivial if $k=n-1$, or $\Pi$ is of type {\rm \textsf{D}},
or complement subsets of distinct types have distinct cardinalities;
however, for some pairs $n,k$ the complement subsets of first and second types have
the same number of elements.}.
Using the characterizations of the collinearity relation in terms
of complement subsets of second type (Lemmas \ref{lemma2-2} and \ref{lemma2-5}),
we establish that $S$ and $U$ are collinear points of ${\mathfrak G}_{k}$
if and only if the same holds for their $f$-images.
In the case (\textsf{C}),
$S$ and $U$ are weak adjacent if and only if $f(S)$ and $f(U)$ are weak adjacent
(this is a consequence of Lemma \ref{lemma-new}).
The required statement follows from Fact \ref{f-2}.

Let $k=0$ and $B=\{p_{1},\dots,p_{2n}\}$ be a base of $\Pi$.
Since for any two points of $\Pi$ there is a base of $\Pi$ containing them,
we need to establish that two points of $B$ are collinear if and only if
their $f$-images are collinear.
As above, we show that $X\subset B$ is a complement subset
if and only if $f(X)$ is a complement subset of the base $f(B)$.
In the case (\textsf{D}), each complement subset is a pair of collinear points
and we get the claim.

Consider the case (\textsf{C}).
For each $i\in \{1,\dots,2n\}$ there exists a point $p$ such that
$$(B\setminus\{p_{i}\})\cup \{p\}\;\mbox{ and }\;(B\setminus\{p_{\sigma(i)}\})\cup \{p\}$$
are bases of $\Pi$ (see proof of Proposition \ref{prop2-1}).
Then
$$(f(B)\setminus\{f(p_{i})\})\cup \{f(p)\}\;\mbox{ and }\;
(f(B)\setminus\{f(p_{\sigma(i)})\})\cup \{f(p)\}$$
are bases of $\Pi$.
Since $f(B)$ is a base of $\Pi$,
there are unique $u,v\in \{1,\dots,2n\}$
such that
$$f(p_{i})\not\perp f(p_{u})\;\mbox{ and }\;f(p_{\sigma(i)})\not\perp f(p_{v}).$$
It is clear that $f(p)$ is non-collinear with $f(p_{u})$ and $f(p_{v})$,
and there is no base of $\Pi$ containing $f(p), f(p_{u}), f(p_{v})$.
This is possible only in the case when $u=\sigma(i), v=i$.
Thus $f(p_{i})\not\perp f(p_{\sigma(i)})$ and $f(p_{i})$ is collinear with
$f(p_{j})$ if $j\ne \sigma(i)$.

\subsection{Inexact and complement subsets of ${\mathcal O}_{\delta}$}
In this subsection we require that $n\ge 4$ and $\Pi$ is of type \textsf{D}.
Let $B=\{p_{1},\dots,p_{2n}\}$ be a base of $\Pi$
and ${\mathcal B}$ be the associated base subset of ${\mathcal O}_{\delta}$
($\delta =+,-$).
Every element of ${\mathcal B}$ contains precisely
one of the points $p_{i}$ or $p_{\sigma(i)}$ for each $i$.

As in subsection 4.1,
we write ${\mathcal B}(+i)$ and ${\mathcal B}(-i)$ for the sets of all
elements of ${\mathcal B}$ which  contain $p_{i}$ or do not contain $p_{i}$, respectively;
for any $i_{1},\dots, i_{s}$ and $j_{1},\dots, j_{u}$
belonging to $\{1,\dots,2n\}$
we define
$${\mathcal B}(+i_{1},\dots, +i_{s},-j_{1},\dots, -j_{u}):=
{\mathcal B}(+i_{1})\cap\dots\cap {\mathcal B}(+i_{s})\cap
{\mathcal B}(-j_{1})\cap \dots\cap{\mathcal B}(-j_{u}).$$
It is clear that
$${\mathcal B}(-i)={\mathcal B}(-i,+\sigma(i))={\mathcal B}(+\sigma(i)).$$

Let ${\mathcal R}\subset {\mathcal B}$.
We say that ${\mathcal R}$ is {\it exact} if
there is only one base subset of ${\mathcal O}_{\delta}$ containing ${\mathcal R}$;
otherwise, ${\mathcal R}$ will be called {\it inexact}.
If ${\mathcal R}\cap {\mathcal B}(+i)$ is not empty then we define
$S_{i}({\mathcal R})$ as the intersection of all subspaces
belonging to ${\mathcal R}$ and containing $p_{i}$,
and we define $S_{i}({\mathcal R}):=\emptyset$ if
the intersection of ${\mathcal R}$ and ${\mathcal B}(+i)$ is empty.
If
$$S_{i}({\mathcal R})=p_{i}$$
for all $i$ then ${\mathcal R}$ is exact;
the converse fails.

\begin{prop}
If $j\ne i,\sigma(i)$ then
\begin{equation}\label{eq3-1}
{\mathcal B}(-i)\cup {\mathcal B}(+i,+j)
\end{equation}
is an inexact subset.
\end{prop}

\begin{proof}
It is clear that
$$p_{j}\in S_{i}({\mathcal B}(-i)\cup {\mathcal B}(+i,+j)).$$
Suppose that $S$ belongs  to \eqref{eq3-1}
and contains $p_{\sigma(j)}$. Then $S$ does not contain $p_{j}$, hence
$p_{i}\not\in S$. The latter means that $p_{\sigma(i)}$ belongs to $S$.
Therefore,
$$p_{\sigma(i)}\in S_{\sigma(j)}({\mathcal B}(-i)\cup {\mathcal B}(+i,+j)).$$
We choose a point $p'_{i}$ on the line $p_{i}p_{j}$
distinct from $p_{i}$ and $p_{j}$.
There is the unique point
$p'_{\sigma(j)}\in p_{\sigma(i)}p_{\sigma(j)}$ collinear with $p'_{i}$.
Then
$$(B\setminus\{p_{i},p_{\sigma(j)}\})\cup\{p'_{i},p'_{\sigma(j)}\}$$
is a base of $\Pi$ and the associated base subset of ${\mathcal O}_{\delta}$
contains \eqref{eq3-1}.
\end{proof}

\begin{prop}\label{prop3-1}
If ${\mathcal R}$ is a maximal inexact subsets of ${\mathcal B}$
then
$${\mathcal R}={\mathcal B}(-i)\cup {\mathcal B}(+i,+j)$$
for some $i$ and $j$ such that $j\ne i,\sigma(i)$.
\end{prop}

\begin{proof}
Since ${\mathcal R}$ is inexact,
we have $S_{i}({\mathcal R})\ne p_{i}$ for certain $i$.
If $S_{i}({\mathcal R})$ is empty  then
${\mathcal R}$ is contained in ${\mathcal B}(-i)$
which contradicts the fact that the inexact subset ${\mathcal R}$
is maximal.
Thus there exists $p_{j}\in S_{i}({\mathcal R})$ such that $j\ne i$.
Then
$${\mathcal R}\subset{\mathcal B}(-i)\cup {\mathcal B}(+i,+j).$$
The inexact subset ${\mathcal R}$ is maximal and the inverse inclusion holds.
\end{proof}

We say that ${\mathcal R}\subset {\mathcal B}$ is a {\it complement subset}
if ${\mathcal B}\setminus {\mathcal R}$
is a maximal inexact subset.
In this case,
$${\mathcal B}\setminus{\mathcal R}={\mathcal B}(-i)\cup {\mathcal B}(+i,+j)$$
(Proposition \ref{prop3-1}) and
$${\mathcal R}={\mathcal B}(+i,+\sigma(j)).$$

\begin{lemma}\label{lemma3-1}
$S,U\in {\mathcal B}$ are collinear points of ${\mathfrak O}_{\delta}$
if and  only if there are precisely
$$\binom{n-2}{2}$$
distinct complement subsets
of ${\mathcal B}$ containing $S$ and $U$.
\end{lemma}

\begin{proof}
See the proof of Lemma \ref{lemma2-2}.
\end{proof}

Let $f$ be a bijective transformation of ${\mathcal O}_{\delta}$
which maps base subsets to base subsets.
Consider $S,U\in {\mathcal O}_{\delta}$ and a base subset
${\mathcal B}\subset{\mathcal O}_{\delta}$ containing them.
As in the previous subsection,
we show that ${\mathcal R}\subset {\mathcal B}$
is a complement subset if and only if $f({\mathcal R})$
is a complement subset $f({\mathcal B})$.
Lemma \ref{lemma3-1} guarantees that
$S$ and $U$ are collinear points of ${\mathfrak O}_{\delta}$
if and only if the same holds for their $f$-images,
and Fact \ref{f-3} gives the claim.

\section{Final remark}
Using the arguments of \cite{Pankov3} (Section 8) and Lemmas \ref{lemma2-1},
\ref{lemma2-2}, \ref{lemma-new}
we can establish the following analogue of the main result of \cite{Pankov3}.
Let $\Pi$ and $\Pi'$ be polar spaces of same rank $n\ge 3$.
We suppose that they both are of type \textsf{C} and
denote by ${\mathcal G}_{k}$ and ${\mathcal G}'_{k}$
the Grassmannians consisting of all $k$-dimensional singular subspaces of $\Pi$ and $\Pi'$,
respectively.
Then every mapping of ${\mathcal G}_{k}$ to ${\mathcal G}'_{k}$
which sends base subsets to base subsets is induced by an embedding of $\Pi$ to $\Pi'$.

\end{document}